\documentclass[12pt,twoside]{article}
\usepackage[english]{babel}
\usepackage[latin1]{inputenc}
\usepackage{amsmath}
\usepackage{amssymb,amsfonts}
\usepackage{graphicx}
\usepackage{times,amssymb,amscd}

\newtheorem{thm}{Theorem}[section]

\newtheorem{conj}[thm]{Conjecture}
\newtheorem{lem}[thm]{Lemma}

\newtheorem{defn}[thm]{Definition}
\newtheorem{rem}[thm]{Remark}
\numberwithin{equation}{section}

\newcommand{\bA}{\mathbf{A}}

\newcommand{\bE}{\mathbf{E}}

\newcommand{\bH}{\mathbf{H}}

\newcommand{\bL}{\mathbf{L}}

\newcommand{\bR}{\mathbf{R}}
\newcommand{\bS}{\mathbf{S}}
\newcommand{\bV}{\mathbf{V}}

\newcommand{\be}{\mathbf{e}}

\newcommand{\br}{\mathbf{r}}
\newcommand{\bx}{\mathbf{x}}
\newcommand{\by}{\mathbf{y}}

\newcommand{\bT}{\mathbf{T}}

\newcommand{\bu}{\mathbf{u}}
\newcommand{\bv}{\mathbf{v}}

\newcommand{\BV}{\boldsymbol{V}}

\newcommand{\Be}{\boldsymbol{e}}
\newcommand{\Bu}{\boldsymbol{u}}
\newcommand{\Bv}{\boldsymbol{v}}

\newcommand{\cP}{\mathcal{P}}
\newcommand{\cS}{\mathcal{S}}

\newcommand{\EUC}{\mathbf E^3}
\newcommand{\SPH}{\bS^3}
\newcommand{\HYP}{\bH^3}
\newcommand{\SXR}{\bS^2\!\times\!\bR}
\newcommand{\HXR}{\bH^2\!\times\!\bR}
\newcommand{\SLR}{\widetilde{\bS\bL_2\bR}}
\newcommand{\NIL}{\mathbf{Nil}}
\newcommand{\SOL}{\mathbf{Sol}}

\begin{document}
\pagestyle{myheadings}
\markboth{\centerline{Jen\H o Szirmai}}
{$\NIL$ geodesic triangles and their interior angle sums}
\title
{$\NIL$ geodesic triangles and \\ their interior angle sums
}

\author{Jen\H o Szirmai \\
\normalsize Budapest University of Technology and \\
\normalsize Economics Institute of Mathematics, \\
\normalsize Department of Geometry \\
\normalsize Budapest, P. O. Box: 91, H-1521 \\
\normalsize szirmai@math.bme.hu
\date{\normalsize{\today}}}

\maketitle
\begin{abstract}

In this paper we study the interior angle sums of geodesic triangles in $\NIL$ geometry
and prove that it can be larger, equal or less than $\pi$.

We use for the computations the projective model of $\NIL$ introduced by E. {Moln\'ar} in \cite{M97}.
\end{abstract}

\section{Introduction} \label{section1}
A geodesic triangle in Riemannian geometry and more generally in metric geometry a
figure consisting of three different points together with the pairwise-connecting geodesic curves.
The points are known as the vertices, while the geodesic curves are known as the sides of the triangle.

In the geometries of constant curvature $\EUC$, $\HYP$, $\SPH$ the well-known sums of the interior angles of geodesic
triangles characterize the space. It is related to the Gauss-Bonnet theorem which states that the integral of the Gauss curvature
on a compact $2$-dimensional Riemannian manifold $M$ is equal to $2\pi\chi(M)$ where $\chi(M)$ denotes the Euler characteristic of $M$.
This theorem has a generalization to any compact even-dimensional Riemannian manifold (see e.g. \cite{Ch}, \cite{KN}).

However, in the other $3$-dimensional homogeneous maximal Riemann spaces (Thur\-ston geo\-metries) there are few results concerning the
angle sums of geodesic triangles. Therefore, it is interesting to study similar question
in the other five Thurston geometries, $\SXR$, $\HXR$, $\NIL$, $\SLR$, $\SOL$.

In \cite{CsSz16} we investigated the angle sum of translation and geodesic triangles in $\SLR$ geometry
and proved that the possible sum of the internal angles in
a translation triangle must be greater or equal than $\pi$. However, in geodesic triangles this sum is
less, greater or equal to $\pi$.

In this paper we consider the analogous problem in $\NIL$ geometry.
\begin{rem}
We note here, that nowadays the $\NIL$ geometry is a widely investigated space concerning
the interesting surfaces, tilings, geodesic and translation ball packings
(see e.g. \cite{B}, \cite{I}, \cite{MSzV}, \cite{MSz06}, \cite{PSchSz11}, \cite{SchSz11}, \cite{Sz07}, \cite{Sz12-1}
and the references given there).
\end{rem}
In \cite{B} K.~Brodaczewska showed, that sum of the interior angles of translation triangles of the $\NIL$ space is larger than $\pi$.

Now, we are interested in {\it geodesic triangles} in $\NIL$ space that is one of the eight Thurston geometries
\cite{S,T}. In Section 2 we describe the projective model of $\NIL$
and we shall use its standard Riemannian metric obtained by
pull back transform to the infinitesimal  arc-length-square at the origin. We also describe the isometry group of $\NIL$,
give an overview about geodesic curves.

In Section 3 we study the $\NIL$ geodesic triangles and prove the interior angle sum of a geodesic triangle
in $\NIL$ geometry can be larger, equal or less than $\pi$.
\section{Basic notions of the $\NIL$ geometry}
In this Section we summarize the significant notions and denotations of the $\NIL$ geometry (see \cite{M97}, \cite{Sz07}).

The $\NIL$ geometry is a homogeneous 3-space derived from the famous real matrix group $\mathbf{L(R)}$ discovered by {W.~Heisenberg}.
The Lie theory with the methode of the projective geometry makes possible to investigate and to describe this topic.

The left (row-column) multiplication of Heisenberg matrices
     \begin{equation}
     \begin{gathered}
     \begin{pmatrix}
         1&x&z \\
         0&1&y \\
         0&0&1 \\
       \end{pmatrix}
       \begin{pmatrix}
         1&a&c \\
         0&1&b \\
         0&0&1 \\
       \end{pmatrix}
       =\begin{pmatrix}
         1&a+x&c+xb+z \\
         0&1&b+y \\
         0&0&1 \\
       \end{pmatrix}
      \end{gathered} \tag{2.1}
     \end{equation}
defines "translations" $\mathbf{L}(\mathbf{R})= \{(x,y,z): x,~y,~z\in \mathbf{R} \}$
on the points of the space $\NIL= \{(a,b,c):a,~b,~c \in \mathbf{R}\}$.
These translations are not commutative in general. The matrices $\mathbf{K}(z) \vartriangleleft \mathbf{L}$ of the form
     \begin{equation}
     \begin{gathered}
       \mathbf{K}(z) \ni
       \begin{pmatrix}
         1&0&z \\
         0&1&0 \\
         0&0&1 \\
       \end{pmatrix}
       \mapsto (0,0,z)
      \end{gathered}\tag{2.2}
     \end{equation}
constitute the one parametric centre, i.e. each of its elements commutes with all elements of $\mathbf{L}$.
The elements of $\mathbf{K}$ are called {\it fibre translations}. $\NIL$ geometry of the Heisenberg group can be projectively
(affinely) interpreted by the "right translations"
on points as the matrix formula
     \begin{equation}
     \begin{gathered}
       (1;a,b,c) \to (1;a,b,c)
       \begin{pmatrix}
         1&x&y&z \\
         0&1&0&0 \\
         0&0&1&x \\
         0&0&0&1 \\
       \end{pmatrix}
       =(1;x+a,y+b,z+bx+c)
      \end{gathered} \tag{2.3}
     \end{equation}
shows, according to (2.1). Here we consider $\mathbf{L}$ as projective collineation
group with right actions in homogeneous coordinates.
We will use the Cartesian homogeneous coordinate simplex $E_0(\be_0)$, $E_1^{\infty}(\be_1)$, $E_2^{\infty}(\be_2)$,
$E_3^{\infty}(\be_3), \ (\{\be_i\}\subset \bV^4$ \ $\text{with the unit point}$ $E(\be = \be_0 + \be_1 + \be_2 + \be_3 ))$
which is distinguished by an origin $E_0$ and by the ideal points of coordinate axes, respectively.
Moreover, $\by=c\bx$ with $0<c\in \mathbf{R}$ (or $c\in\mathbf{R}\setminus\{0\})$
defines a point $(\bx)=(\by)$ of the projective 3-sphere $\cP \cS^3$ (or that of the projective space $\cP^3$ where opposite rays
$(\bx)$ and $(-\bx)$ are identified).
The dual system $\{(\Be^i)\}, \ (\{\Be^i\}\subset \BV_4)$ describes the simplex planes, especially the plane at infinity
$(\Be^0)=E_1^{\infty}E_2^{\infty}E_3^{\infty}$, and generally, $\Bv=\Bu\frac{1}{c}$ defines a plane
$(\Bu)=(\Bv)$ of $\cP \cS^3$
(or that of $\cP^3$). Thus $0=\bx\Bu=\by\Bv$ defines the incidence
of point $(\bx)=(\by)$ and plane
$(\Bu)=(\Bv)$.
Thus {\bf Nil} can be visualized in the affine 3-space $\bA^3$
(so in $\bE^3$) as well.

In this context E. Moln\'ar \cite{M97} has derived the well-known infinitesimal arc-length square invariant under translations $\bL$ at any point of $\NIL$ as follows
\begin{equation}
   \begin{gathered}
      (dx)^2+(dy)^2+(-xdy+dz)^2=\\
      (dx)^2+(1+x^2)(dy)^2-2x(dy)(dz)+(dz)^2=:(ds)^2
       \end{gathered} \tag{2.4}
     \end{equation}
Hence we get the symmetric metric tensor field $g$ on $\NIL$ by components $g_{ij}$, furthermore its inverse:
\begin{equation}
   \begin{gathered}
       g_{ij}:=
       \begin{pmatrix}
         1&0&0 \\
         0&1+x^2&-x \\
         0&-x&1 \\
         \end{pmatrix},  \quad  g^{ij}:=
       \begin{pmatrix}
         1&0&0 \\
         0&1&x \\
         0&x&1+x^2 \\
         \end{pmatrix} \\
         \text{with} \ \det(g_{ij})=1.
        \end{gathered} \tag{2.5}
     \end{equation}
The translation group $\mathbf{L}$ defined by formula (2.3) can be extended to a larger group $\mathbf{G}$ of collineations,
preserving the fibering, that will be equivalent to the (orientation preserving) isometry group of $\NIL$.
In \cite{M06} E.~Moln\'ar has shown that
a rotation trough angle $\omega$
about the $z$-axis at the origin, as isometry of $\NIL$, keeping invariant the Riemann
metric everywhere, will be a quadratic mapping in $x,y$ to $z$-image $\overline{z}$ as follows:
     \begin{equation}
     \begin{gathered}
       \br(O,\omega):(1;x,y,z) \to (1;\overline{x},\overline{y},\overline{z}); \\
       \overline{x}=x\cos{\omega}-y\sin{\omega}, \ \ \overline{y}=x\sin{\omega}+y\cos{\omega}, \\
       \overline{z}=z-\frac{1}{2}xy+\frac{1}{4}(x^2-y^2)\sin{2\omega}+\frac{1}{2}xy\cos{2\omega}.
      \end{gathered} \tag{2.6}
     \end{equation}
This rotation formula, however, is conjugate by the quadratic mapping $\mathcal{M}$
     \begin{equation}
     \begin{gathered}
       x \to x'=x, \ \ y \to y'=y, \ \ z \to z'=z-\frac{1}{2}xy  \ \ \text{to} \\
       (1;x',y',z') \to (1;x',y',z')
       \begin{pmatrix}
         1&0&0&0 \\
         0&\cos{\omega}&\sin{\omega}&0 \\
         0&-\sin{\omega}&\cos{\omega}&0 \\
         0&0&0&1 \\
       \end{pmatrix}
       =(1;x",y",z"), \\
       \text{with} \ \ x" \to \overline{x}=x", \ \ y" \to \overline{y}=y", \ \ z" \to \overline{z}=z"+\frac{1}{2}x"y",
      \end{gathered} \tag{2.7}
     \end{equation}
i.e. to the linear rotation formula. This quadratic conjugacy modifies the $\NIL$ translations in (2.3), as well.
We shall use the following important classification theorem.
\begin{thm}[E.~Moln\'ar \cite{M06}]
\begin{enumerate}
\item Any group of $\NIL$ isometries, containing a 3-dimensional translation lattice,
is conjugate by the quadratic mapping in (2.5) to an affine group of the affine (or Euclidean) space $\bA^3=\EUC$
whose projection onto the (x,y) plane is an isometry group of $\bE^2$. Such an affine group preserves a plane
$\to$ point polarity of signature $(0,0,\pm 0,+)$.
\item Of course, the involutive line reflection about the $y$ axis
     \begin{equation}
     \begin{gathered}
       (1;x,y,z) \to (1;-x,y,-z),
      \end{gathered} \notag
     \end{equation}
preserving the Riemann metric, and its conjugates by the above isometries in {$1$} (those of the identity component)
are also {$\NIL$}-isometries. There does not exist orientation reversing $\NIL$-isometry.
\end{enumerate}
\end{thm}
\begin{rem}
We obtain from the above described projective model a new model of $\NIL$ geometry derived by the quadratic mapping $\mathcal{M}$.
This is the {\it linearized model of $\NIL$ space} (see \cite{B}).
\end{rem}
\subsection{Geodesic curves} \label{subsection2}
The geodesic curves of the $\NIL$ geometry are generally defined as having locally minimal arc length between their any two (near enough) points.
The equation systems of the parametrized geodesic curves $g(x(t),y(t),z(t))$  in our model can be determined by the
general theory of Riemann geometry.
We can assume, that the starting point of a geodesic curve is the origin because we can transform a curve into an
arbitrary starting point by translation (2.1);
\begin{equation}
\begin{gathered}
        x(0)=y(0)=z(0)=0; \ \ \dot{x}(0)=c \cos{\alpha}, \ \dot{y}(0)=c \sin{\alpha}, \\ \dot{z}(0)=w; \ - \pi \leq \alpha \leq \pi. \notag
\end{gathered}
\end{equation}
The arc length parameter $s$ is introduced by
\begin{equation}
 s=\sqrt{c^2+w^2} \cdot t, \ \text{where} \ w=\sin{\theta}, \ c=\cos{\theta}, \ -\frac{\pi}{2}\le \theta \le \frac{\pi}{2}, \notag
\end{equation}
i.e. unit velocity can be assumed.

The equation systems of a helix-like geodesic curves $g(x(t),y(t),z(t))$ if $0<|w| <1 $:
\begin{equation}
\begin{gathered}
x(t)=\frac{2c}{w} \sin{\frac{wt}{2}}\cos\Big( \frac{wt}{2}+\alpha \Big),\ \
y(t)=\frac{2c}{w} \sin{\frac{wt}{2}}\sin\Big( \frac{wt}{2}+\alpha \Big), \notag \\
z(t)=wt\cdot \Big\{1+\frac{c^2}{2w^2} \Big[ \Big(1-\frac{\sin(2wt+2\alpha)-\sin{2\alpha}}{2wt}\Big)+ \\
+\Big(1-\frac{\sin(2wt)}{wt}\Big)-\Big(1-\frac{\sin(wt+2\alpha)-\sin{2\alpha}}{2wt}\Big)\Big]\Big\} = \\
=wt\cdot \Big\{1+\frac{c^2}{2w^2} \Big[ \Big(1-\frac{\sin(wt)}{wt}\Big)
+\Big(\frac{1-\cos(2wt)}{wt}\Big) \sin(wt+2\alpha)\Big]\Big\}. \tag{2.8}
\end{gathered}
\end{equation}
In the cases $w=0$ the geodesic curve is the following:
\begin{equation}
x(t)=c\cdot t \cos{\alpha}, \ \ y(t)=c\cdot t \sin{\alpha}, \ \ z(t)=\frac{1}{2} ~ c^2 \cdot t^2 \cos{\alpha} \sin{\alpha}. \tag{2.9}
\end{equation}
The cases $|w|=1$ are trivial: $(x,y)=(0,0), \ z=w \cdot t$.
\begin{defn}
The distance $d(P_1,P_2)$ between the points $P_1$ and $P_2$ is defined by the arc length of geodesic curve
from $P_1$ to $P_2$.
\end{defn}
\section{Geodesic triangles} \label{section3}
We consider $3$ points $A_1$, $A_2$, $A_3$ in the projective model of $\NIL$ space (see Section 2).
The {\it geodesic segments} $a_k$ between the points $A_i$ and $A_j$
$(i<j,~i,j,k \in \{1,2,3\}, k \ne i,j$) are called sides of the {\it geodesic triangle} with vertices $A_1$, $A_2$, $A_3$.
\begin{figure}[ht]
\centering
\includegraphics[width=12cm]{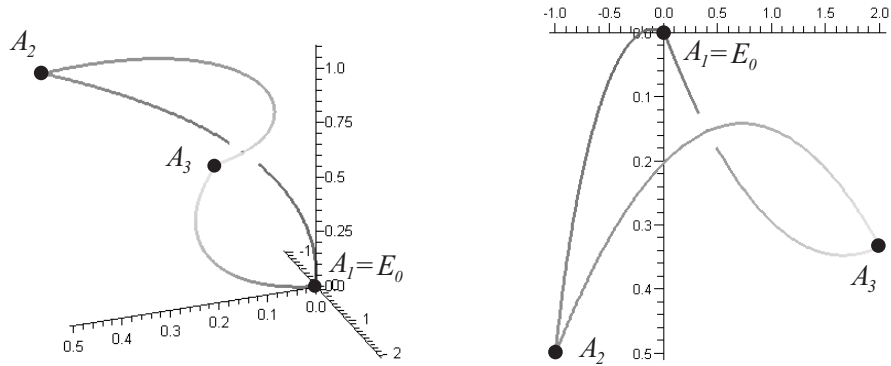}
\caption{Geodesic triangle with vertices $A_1=(1,0,0,0)$, $A_2=(1,1/2,-1,1)$, $A_3=(1,1/3,2,1)$ where its interior angle sum
is $\approx 3.45294 > \pi $}.
\label{}
\end{figure}
In Riemannian geometries the metric tensor (see (2.5)) is used to define the angle $\theta$ between two geodesic curves.
If their tangent vectors in their common point are $\bu$ and $\bv$ and $g_{ij}$ are the components of the metric tensor then
\begin{equation}
\cos(\theta)=\frac{u^i g_{ij} v^j}{\sqrt{u^i g_{ij} u^j~ v^i g_{ij} v^j}} \tag{3.1}
\end{equation}
It is clear by the above definition of the angles and by the metric tensor (2.5), that
the angles are the same as the Euclidean ones at the origin of
the projective model of $\NIL$ geometry.

We note here that the angle of two intersecting geodesic curves depend on the orientation of the tangent vectors. We will consider
the {\it internal angles} of the triangles that are denoted at the vertex $A_i$ by $\omega_i$ $(i\in\{1,2,3\})$.

\subsection{Fibre-like right angled triangles}
A geodesic triangle is called fibre-like if one of its edges lies on a fibre line. In this section we study the right angled fibre-like
triangles. We can assume without loss of generality
that the veritices $A_1$, $A_2$, $A_3$ of a fibre-like right angled triangle (see Fig.~2) have the following coordinates:
\begin{equation}
A_1=(1,0,0,0),~A_2=(1,0,0,z^2),~A_3=(1,x^3,0,z^3=z^2) \tag{3.2}
\end{equation}
\begin{figure}[ht]
\centering
\includegraphics[width=13cm]{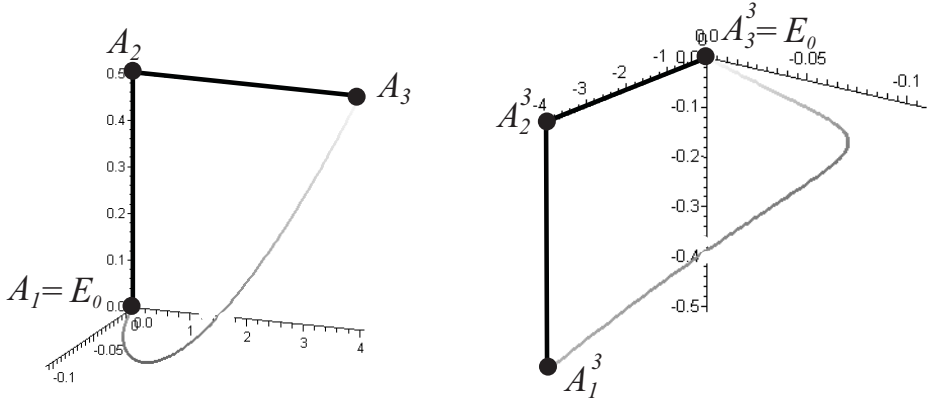}
\caption{Fibre-like geodesic triangle $A_1A_2A_3$ and its translated image
$A_1^3A_2^3A_3^3$ by translation $\bT_{A_3}$ where $A_1=(1,0,0,0)$,~$A_2=(1,0,0,\frac{1}{2})$,~$A_3=(1,4,0,\frac{1}{2})$.}
\label{}
\end{figure}
The geodesic segment $A_2A_3$ lies on a to $x$ axis parallel straight line, the geodesic segment $A_1A_2$ lies on the $z$ axis
and their angle is $\omega_2=\frac{\pi}{2}$ in the $\NIL$ space (this angle is in Euclidean sense also $\frac{\pi}{2}$) (see Fig.~2).

{\it In order to determine the further internal angles of the fibre-like geodesic triangle} $A_1A_2A_3$
we define \emph{translations} $\bT_{A_i}$, $(i\in \{2,3\})$ as elements of the isometry group of $\NIL$, that
maps the origin $E_0$ onto $A_i$ (see Fig.~2).
E.g. the isometrie $\bT_{A_3}$ and its inverse (up to a positive determinant factor) can be given by:
\begin{equation}
\bT_{A_3}=
\begin{pmatrix}
1 & x^3 & 0 & z^3 \\
0 & 1 & 0 & 0 \\
0 & 0 & 1 & x^3 \\
0 & 0 & 0 & 1
\end{pmatrix} , ~ ~ ~
\bT_{A_3}^{-1}=
\begin{pmatrix}
1 & -x^3 & 0 & -z^3 \\
0 & 1 & 0 & 0 \\
0 & 0 & 1 & -x^3 \\
0 & 0 & 0 & 1
\end{pmatrix} , \tag{3.3}
\end{equation}
and the images $\bT_{A_3}(A_i)$ of the vertices $A_i$ $(i \in \{1,2,3\})$ are the following (see also Fig.~2):
\begin{equation}
\begin{gathered}
\bT^{-1}_{A_3}(A_1)=A_1^3=(1,-x^3,0,-z^3),~\bT^{-1}_{A_3}(A_2)=A_2^3=(1,-x^3,0,0), \\ \bT^{-1}_{A_3}(A_3)=A_3^3=E_0=(1,0,0,0). \tag{3.4}
\end{gathered}
\end{equation}
\begin{figure}[ht]
\centering
\includegraphics[width=12cm]{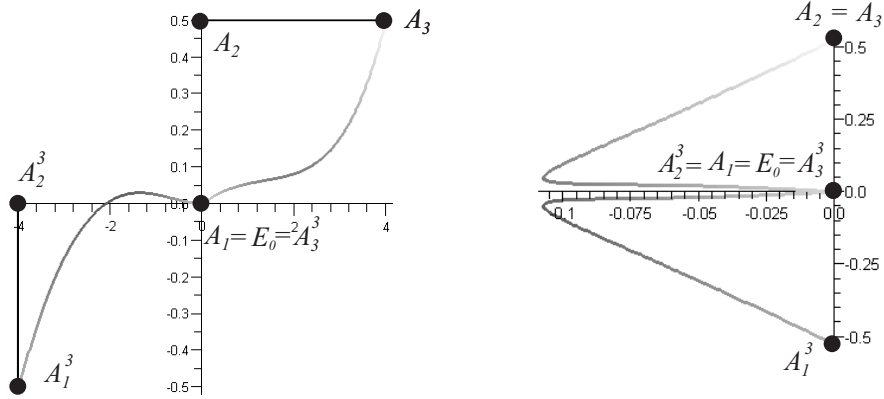}
\caption{The points $A_3$ and $A_1^3$ are antipodal related to the origin $E_0$
($A_1=(1,0,0,0)$,~$A_2=(1,0,0,\frac{1}{2})$,~$A_3=(1,4,0,\frac{1}{2})$.)}
\label{}
\end{figure}
Our aim is to determine angle sum $\sum_{i=1}^3(\omega_i)$ of the internal angles of the above right angled
fibre-like geodesic triangle $A_1A_2A_3$.
We have seen that $\omega_2=\frac{\pi}{2}$ and the angle of geodesic curves with common point at the origin $E_0$ is the same as the
Euclidean one therefore it can be determined by usual Euclidean sense. Moreover, the translation $\bT_{A_3}$ is isometry
in $\NIL$ geometry thus
$\omega_3$ is equal to the angle $(g(A_3^3, A_1^3)g(A_3^3, A_2^3)) \angle$ (see Fig.~2,3)
where $g(A_3^3, A_1^3)$, $g(A_3^3, A_2^3)$ are oriented geodesic curves $(E_0=A_2^2=A_3^3)$ and
$\omega_1$ is equal to the angle $(g(E_0, A_3)g(E_0, A_2)) \angle$ (see Fig.~2,3)
where $g(E_0, A_3)$, $g(E_0, A_2)$ are also oriented geodesic curves.

We denote the oriented unit tangent vectors of the oriented geodesic curves $g(E_0, A_i^j)$ with $\mathbf{t}_i^j$ where
$(i,j)=(1,3),(2,3),(3,0),(2,0)$ and $A_3^0=A_3$, $A_2^0=A_2$.
The Euclidean coordinates of $\mathbf{t}_i^j$ (see Section 2.1) are :
\begin{equation}
\mathbf{t}_i^j=(\cos(\theta_i^j) \cos(\alpha_i^j), \cos(\theta_i^j) \sin(\alpha_i^j), \sin(\theta_i^j)). \tag{3.5}
\end{equation}
\begin{lem}
The sum of the interior angles of a fibre-like right angled geodesic triangle is greather or equal to $\pi$.
\end{lem}
\textbf{Proof:} It is clear, that $\mathbf{t}_2^0=(0,0,1)$ and $\mathbf{t}_2^3=(-1,0,0)$. Moreover,
the points $A_3$ and $A_1^3$ are antipodal related to the origin $E_0$ therefore the equation $|\theta_3^0|=|\theta_1^3|$ holds
(i.e. the angle betveen the vector $\mathbf{t}_3^0$ and $[x,y]$ plane are equal to the angle between the vector $\mathbf{t}_1^3$ and the $[x,y]$
plane). Moreover, we have seen, that $\omega_2=\frac{\pi}{2}$.
That means, that $\omega_1=\frac{\pi}{2}-|\theta_3^0|=\frac{\pi}{2}-|\theta_1^3|$.

The vector $\mathbf{t}_2^3$ lies in the $[x,y]$ plane therefore the angle $\omega_3$ is greather or equal than $|\theta_3^0|=|\theta_1^3|$.
Finally we obtain, that
$$
\sum_{i=1}^3(\omega_i)=\frac{\pi}{2}+\frac{\pi}{2}-|\theta_3^0|+\omega_3 \ge \pi. \ \ \square
$$
\vspace{3mm}
We fix the $z^2=z^3\in \mathbf{R}$ coordinates of the vertices $A_2$ and $A_3$ and study the the internal angle sum
$\sum_{i=1}^3(\omega_i(x^3))$ of the right angled
geodesic triangle $A_1A_2A_3$ if $x^3$ coordinate of vertex $A_3$ tends to zero or infinity.
We obtain directly from the system of equation (2.8) of geodesic curves the following
\begin{lem}
If the coordinates $z^2=z^3 \in \mathbf{R}$ are fixed then
$$\lim_{x^3 \to 0}(\omega_1(x^3))=0, ~
\lim_{x^3 \to 0}(\omega_3(x^3))=\frac{\pi}{2} ~ \Rightarrow ~ \lim_{x^3 \to 0} \Bigg(\sum_{i=1}^3(\omega_i(x^3))\Bigg)=\pi,$$
$$\lim_{x^3 \to \infty}(\omega_1(x^3))=\frac{\pi}{2}, ~
\lim_{x^3 \to \infty}(\omega_3(x^3))=0 ~ \Rightarrow ~ \lim_{x^3 \to \infty} \Bigg(\sum_{i=1}^3(\omega_i(x^3))\Bigg)=\pi.$$
\end{lem}
In the following table we summarize some numerical data of geodesic triangles for given parameters:
\medbreak
\centerline{\vbox{
\halign{\strut\vrule\quad \hfil $#$ \hfil\quad\vrule
&\quad \hfil $#$ \hfil\quad\vrule &\quad \hfil $#$ \hfil\quad\vrule &\quad \hfil $#$ \hfil\quad\vrule &\quad \hfil $#$ \hfil\quad\vrule &\quad \hfil $#$ \hfil\quad\vrule
\cr
\noalign{\hrule}
\multispan6{\strut\vrule\hfill\bf Table 1, ~ ~  $z^2=z^3=1/2$ \hfill\vrule}%
\cr
\noalign{\hrule}
\noalign{\vskip2pt}
\noalign{\hrule}
x^3 & |\theta_3^0|=|\theta_1^3| & d(A_1A_3) & \omega_1 & \omega_3  & \sum_{i=1}^3(\omega_i)  \cr
\noalign{\hrule}
\rightarrow 0 & \rightarrow \pi/2  & 1/2 & \rightarrow 0 & \rightarrow \pi/2 & \rightarrow \pi \cr
\noalign{\hrule}
1/1000 & 1.56878 & 0.50000 & 0.00202 & 1.56884 & 3.14166 \cr
\noalign{\hrule}
1/3 & 0.97374 & 0.59901 & 0.59705 & 0.99435 & 3.16220 \cr
\noalign{\hrule}
1 & 0.42883 & 1.10937 & 1.14197 & 0.48351 & 3.19627 \cr
\noalign{\hrule}
4 & 0.05337 & 4.01337 & 1.51743 & 0.11957 & 3.20779 \cr
\noalign{\hrule}
15 & 0.00169 & 15.00042 & 1.56911 & 0.01277 & 3.15268 \cr
\noalign{\hrule}
100 & 0.00001 & 100.00000 & 1.57079 & 0.00030 & 3.14189 \cr
\noalign{\hrule}
\rightarrow \infty & \rightarrow 0 & \rightarrow \infty & \rightarrow \pi/2 & \rightarrow 0 & \rightarrow \pi \cr
\noalign{\hrule}
}}}
\medbreak
\subsubsection{Hyperbolic-like right angled geodesic triangles}
A geodesic triangle is hyperbolic-like if its vertices lie in the base plane (i.e. $[x,y]$ coordinate plane) of the model.
In this section we analyse the internal angle sum of the right angled hyperbolic-like
triangles. We can assume without loss of generality
that the veritices $A_1$, $A_2$, $A_3$ of a hyperbolic-like right angled triangle (see Fig.~4) $T_g$ have the following coordinates:
\begin{equation}
A_1=(1,0,0,0),~A_2=(1,0,y^2,0),~A_3=(1,x^3,y^2=y^3,0) \tag{3.6}
\end{equation}
\begin{figure}[ht]
\centering
\includegraphics[width=12cm]{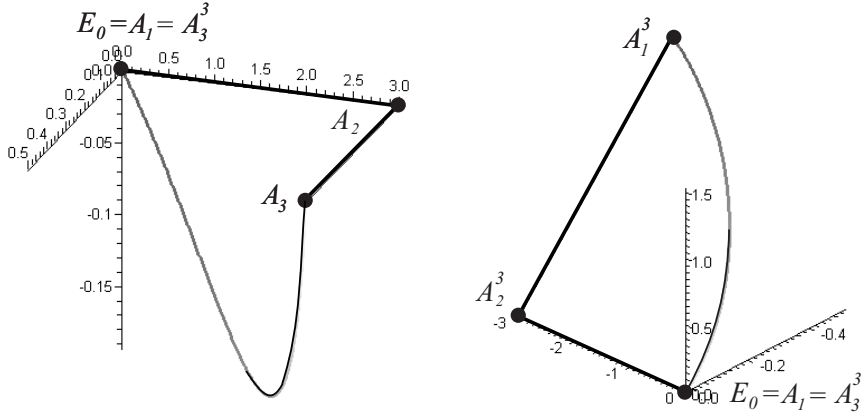}
\caption{Hyperbolic-like geodesic triangle $A_1A_2A_3$ and its tanslated copy $A_1^3A_2^3A_3^3$.}
\label{}
\end{figure}
The geodesic segment $A_1A_2$ lies on the $y$ axis, the geodesic segment $A_2A_3$ lies parallel to the $x$ axis containing the point $A_2$.
It is clear that $\omega_2=\frac{\pi}{2}$ in the $\NIL$ space (this angle is in Euclidean sense also $\frac{\pi}{2}$).

{\it In order to determine the further internal angles of the fibre-like geodesic triangle} $A_1A_2A_3$
similarly to the fibre-like case we define the \emph{translation} $\bT_{A_3}$, (see (2.6)) that
maps the origin $E_0$ onto $A_3$ that can be given by:
\begin{equation}
\bT_{A_3}=
\begin{pmatrix}
1 & x^3 & y^2 & 0 \\
0 & 1 & 0 & 0 \\
0 & 0 & 1 & x^3 \\
0 & 0 & 0 & 1
\end{pmatrix} , ~ ~ ~
\bT_{A_3}^{-1}=
\begin{pmatrix}
1 & -x^3 & -y^2 & x^3y^2 \\
0 & 1 & 0 & 0 \\
0 & 0 & 1 & -x^3 \\
0 & 0 & 0 & 1
\end{pmatrix}. \tag{3.7}
\end{equation}
We obtain that the images $\bT^{-1}_{A_3}(A_i)$ of the vertex $A_i$ $(i \in \{1,2,3\}),$
are the following (see also Fig.~4):
\begin{equation}
\begin{gathered}
\bT^{-1}_{A_3}(A_1)=A_1^3=(1,-x^3,-y^2,x^3y^2),~\bT^{-1}_{A_3}(A_2)=A_2^3=(1,-x^3,0,0), \\ \bT^{-1}_{A_3}(A_3)=A_3^3=E_0=(1,0,0,0). \tag{3.8}
\end{gathered}
\end{equation}
We study similarly to the above fibre-like case the sum $\sum_{i=1}^3(\omega_i)$ of the interior angles of the above right angled
hyperbolic-like geodesic triangle $A_1A_2A_3$.

It is clear, that the angle of geodesic curves with common point at the origin $E_0$ is the same as the
Euclidean one therefore it can be determined by usual Euclidean sense. The translation $\bT_{A_3}$ preserve
the measure of angles $\omega_i$ $(i \in \{2,3\})$ therefore (see Fig.~4)
$\omega_3=(g(A_3^3 A_1^3)g(A_3^3 A_2^3)) \angle$ ($A_2^2=A_3^3=E_0$ and $g(E_0A_i^j)$ $((i,j)=(1,3), (2,3))$ are oriented geodesic curves).

Similarly to the fibre-like case the Euclidean coordinates of the oriented unit tangent vector $\mathbf{t}_i^j$ of the oriented geodesic curves
$g(E_0, A_i^j)$ $((i,j)=(2,3),(3,2),$ $(1,2),(1,3))$ is given by (3.5).

First we fix the $x^3\in \mathbf{R}$ coordinate of the vertices $A_3$ and study the the internal angle sum
$\sum_{i=1}^3(\omega_i(y^2=y^3))$ of the right angled
geodesic triangle $A_1A_2A_3$ if $y^2=y^3$ coordinates of vertices $A_2$ and $A_3$ tend to zero or infinity.
We obtain directly from the system of equation (2.8) of geodesic curves the following
\begin{lem}
If the coordinate $x^3 \in \mathbf{R}$ is fixed then
$$\lim_{y^2=y^3 \to 0}(\omega_1(y^2))=\frac{\pi}{2}, ~
\lim_{y^2=y^3 \to 0}(\omega_3(y^2))=0 ~ \Rightarrow ~ \lim_{y^2=y^3 \to 0} \Bigg(\sum_{i=1}^3(\omega_i(y^2))\Bigg)=\pi,$$
$$\lim_{y^2=y^3 \to \infty}(\omega_1(y^2))=0, ~
\lim_{y^2=y^3 \to \infty}(\omega_3(y^2))=\frac{\pi}{2}~ \Rightarrow ~ \lim_{y^2=y^3 \to \infty} \Bigg(\sum_{i=1}^3(\omega_i(y^2))\Bigg)=\pi.$$
\end{lem}
Secondly we fix the $y^2=y^3\in \mathbf{R}$ coordinates of the vertices $A_2$ and $A_3$ and study the internal angle sum
$\sum_{i=1}^3(\omega_i(x^3))$ of the right angled
geodesic triangle $A_1A_2A_3$ if $x^3$ coordinate of vertex $A_3$ tends to zero or infinity.
We obtain directly from the system of equation (2.8) of geodesic curves the following
\begin{lem}
If the coordinates $y^2=y^3 \in \mathbf{R}$ are fixed then
$$\lim_{x^3 \to 0}(\omega_1(x^3))=0, ~
\lim_{x^3 \to 0}(\omega_3(x^3))=\frac{\pi}{2} ~\Rightarrow ~ \lim_{x^3 \to 0} \Bigg(\sum_{i=1}^3(\omega_i(x^3))\Bigg)=\pi,$$
$$\lim_{x^3 \to \infty}(\omega_1(x^3))=\frac{\pi}{2}, ~
\lim_{x^3 \to \infty}(\omega_3(x^3))=0 ~ \Rightarrow ~ \lim_{x^3 \to \infty} \Bigg(\sum_{i=1}^3(\omega_i(x^3))\Bigg)=\pi.$$
\end{lem}
We can determine the interior angle sum of arbitrary hyperbolic-like geodesic triangle similarly to the fibre-like case.
In the following table we summarize some numerical data of geodesic triangles for given parameters:
\medbreak
\centerline{\vbox{
\halign{\strut\vrule\quad \hfil $#$ \hfil\quad\vrule
&\quad \hfil $#$ \hfil\quad\vrule &\quad \hfil $#$ \hfil\quad\vrule &\quad \hfil $#$ \hfil\quad\vrule &\quad \hfil $#$ \hfil\quad\vrule &\quad \hfil $#$ \hfil\quad\vrule
\cr
\noalign{\hrule}
\multispan6{\strut\vrule\hfill\bf Table 2, ~ ~  $x^3=1/2$ \hfill\vrule}%
\cr
\noalign{\hrule}
\noalign{\vskip2pt}
\noalign{\hrule}
y^2=y^3 & |\theta_3^0|=|\theta_1^3| & d(A_1A_3) & \omega_1 & \omega_3  & \sum_{i=1}^3(\omega_i)  \cr
\noalign{\hrule}
\rightarrow 0 & \rightarrow \pi/2  & 1/2 & \rightarrow \pi/2 & \rightarrow 0 & \rightarrow \pi \cr
\noalign{\hrule}
1/100 & 0.00490 & 0.50011 & 1.54958 & 0.01940 & 3.13977 \cr
\noalign{\hrule}
1/3 & 0.13378 & 0.60651 & 0.94882 & 0.56204 & 3.08166 \cr
\noalign{\hrule}
3 & 0.13700 & 3.09310 & 0.14449 & 1.19813 & 2.91342 \cr
\noalign{\hrule}
6 & 0.06132 & 6.06701 & 0.11959 & 1.30229 & 2.99268 \cr
\noalign{\hrule}
20 & 0.00726 & 20.02442 & 0.04828 & 1.47308 & 3.09216 \cr
\noalign{\hrule}
100 & 0.00030 & 100.00500 & 0.00999 & 1.55082 & 3.13160 \cr
\noalign{\hrule}
\rightarrow \infty & \rightarrow 0 & \rightarrow \infty & \rightarrow 0 & \rightarrow \pi/2 & \rightarrow \pi \cr
\noalign{\hrule}
}}}
\medbreak
\medbreak
\centerline{\vbox{
\halign{\strut\vrule\quad \hfil $#$ \hfil\quad\vrule
&\quad \hfil $#$ \hfil\quad\vrule &\quad \hfil $#$ \hfil\quad\vrule &\quad \hfil $#$ \hfil\quad\vrule &\quad \hfil $#$ \hfil\quad\vrule &\quad \hfil $#$ \hfil\quad\vrule
\cr
\noalign{\hrule}
\multispan6{\strut\vrule\hfill\bf Table 3, ~ ~  $y^2=y^3=1/3$ \hfill\vrule}%
\cr
\noalign{\hrule}
\noalign{\vskip2pt}
\noalign{\hrule}
x^3 & |\theta_3^0|=|\theta_1^3| & d(A_1A_3) & \omega_1 & \omega_3  & \sum_{i=1}^3(\omega_i)  \cr
\noalign{\hrule}
\rightarrow 0 & \rightarrow \pi/2  & 1/3 & \rightarrow 0 & \rightarrow \pi/2 & \rightarrow \pi \cr
\noalign{\hrule}
1/100 & 0.00495 & 0.33349 & 0.02958 & 1.53998 & 3.14036 \cr
\noalign{\hrule}
1/3 & 0.11518 & 0.47461 &  0.76510 & 0.76510 & 3.10100  \cr
\noalign{\hrule}
3 & 0.09346 & 3.04190 & 1.31933 & 0.09854 & 2.98866 \cr
\noalign{\hrule}
6 & 0.04131 & 6.02995 & 1.39094 & 0.08042 & 3.04215 \cr
\noalign{\hrule}
20 & 0.00485 & 20.01086 & 1.50562 & 0.03222 & 3.10863 \cr
\noalign{\hrule}
100 & 0.00020 & 100.00222 & 1.55748 & 0.00666 & 3.13493 \cr
\noalign{\hrule}
\rightarrow \infty & \rightarrow 0 & \rightarrow \infty & \rightarrow \pi/2 & \rightarrow 0 & \rightarrow \pi \cr
\noalign{\hrule}
}}}
\medbreak
Finally, we get the following Lemma
\begin{lem}
The interior angle sums of hyperbolic-like geodesic triangles can be less or equal to $\pi$.
\end{lem}
\begin{conj}
The sum of the interior angles of any hyperbolic-like right angled hyperbolic-like right angled geodesic triangle is less or equal to $\pi$.
\end{conj}
\subsubsection{Geodesic triangles with internal angle sum $\pi$}
In the above sections we discuss the fibre- and hyperbolic-like geodesic triangles and proved that there are right angled
geodesic triangles whose angle sum $\sum_{i=1}^{3}(\omega_i)$ is greather, less or equal to $\pi$, but $\sum_{i=1}^{3}(\omega_i)=\pi$ is
realized if one of the vertices of a geodesic triangle $A_1A_2A_3$ tends to the infinity (see Tables 1,2,3). We prove the following
\begin{lem}
There is geodesic triangle $A_1A_2A_3$ with internal angle sum $\pi$ where its vertices are {\it proper}
(i.e. $A_i \in \NIL$, ($i=1,2,3$)).
\end{lem}
{\bf{Proof:}} We consider a hyperbolic-like geodesic right angled triangle with vertices $A_1=E_0=(1,0,0,0)$, $A_2^h=(1,0,y^2_h,0)$
$A_3^h=(1,x^3_h,y^2_h=y^3_h,0)$ and a fibre-like right angled geodesic triangle with vertices
$A_1=E_0=(1,0,0,0)$, $A_2^f=(1,0,0,z^2_f=z^3_f)$ $A_3^f=(1,x^3_f,0,z^2_f=z^3_f)$
($0 < x^3_h, x^3_f, y^2_h=y^3_h, z^2_f=z^3_f < \infty $). We cosider the straight line segment (in Euclidean sense) $A_3^f A_3^h \subset \NIL$.

We consider a geodesic right angled triangle $A_1A_2A_3(t)$ where $A_3(t) \in A_3^f A_3^h$, $(t\in [0,1])$. $A_3(t)$ is moving on the
segment $A_3^hA_3^f$ and if $t=0$ then $A_3(0)=A_3^h$, if $t=1$ then $A_3(1)=A_3^f$.

Similarly to the above cases the internal angles of the geodesic triangle $A_1A_2$ $A_3(t)$ are denoted by $\omega_i(t)$ $(i \in \{1,2,3\})$.
The angle sum $\sum_{i=1}^{3}(\omega_i(0)) < \pi$ and $\sum_{i=1}^{3}(\omega_i(1)) > \pi$. Moreover the angles $\omega_i(t)$ change
continuously if the parameter $t$ run in the interval $[0,1]$. Therefore there is a $t_E\in (0,1)$ where $\sum_{i=1}^{3}(\omega_i(t_E)) = \pi$.
~ ~ $\square$

We obtain by the Lemmas of this Section the following
\begin{thm}
The sum of the internal angles of a geodesic triangle of $\NIL$ space can be greather, less or equal to $\pi$.
\end{thm}
%

\end{document}